\documentclass[12pt,a4paper]{amsart}
\usepackage[utf8]{inputenc}	
\usepackage[T1]{fontenc}
\usepackage[in,headings]{fullpage}
\usepackage{amsbsy, amscd, amsfonts, amsmath, amsrefs, amssymb, amsthm}
\usepackage{graphicx}
\usepackage{float}
\usepackage{color}   
\usepackage[colorlinks=true,linkcolor=blue,citecolor=blue,urlcolor=blue]{hyperref}

\newtheorem{theorem}{Theorem}

\newtheorem{lemma}[theorem]{Lemma}

\theoremstyle{definition}

\theoremstyle{remark}

\newcommand{\R}{\mathbf{R}}

\renewcommand{\Re}{\mathop{\mathrm{Re}}\nolimits}

\newcommand{\Rzeta}{\mathop{\mathcal R }\nolimits}

\newfont{\cmbsy}{cmbsy10}
\newfont{\cmmib}{cmmib10}
\newcommand{\Orden}{\mathop{\hbox{\cmbsy O}}\nolimits}
\newcommand{\orden}{\mathop{\hbox{\cmmib o}}\nolimits}

\begin{document}

\title[Mean Values of the auxiliary function]
{Mean Values of the auxiliary function.}
\author[Arias de Reyna]{J. Arias de Reyna}
\address{%
Universidad de Sevilla \\ 
Facultad de Matem\'aticas \\ 
c/Tarfia, sn \\ 
41012-Sevilla \\ 
Spain.} 




\email{arias@us.es, ariasdereyna1947@gmail.com}


\begin{abstract} Let $\Rzeta(s)$ be the function related to $\zeta(s)$ found by Siegel in the papers of Riemann. In this paper we obtain the main terms of the mean values \[\frac{1}{T}\int_0^T |\Rzeta(\sigma+it)|^2\Bigl(\frac{t}{2\pi}\Bigr)^\sigma\,dt, \quad\text{and}\quad \frac{1}{T}\int_0^T |\Rzeta(\sigma+it)|^2\,dt.\] Giving complete proofs of some result of the paper of Siegel about the Riemann Nachlass.  Siegel follows Riemann to obtain these mean values. We have followed a more standard path, and explain the difficulties we encountered in understanding Siegel's reasoning.

\end{abstract}

\maketitle

\section{Introduction} 

This paper with  \cite{A98}, \cite{A100}, \cite{A193} and \cite{A102} attempts to give complete proofs of the developments and results on the zeros of the auxiliary function contained in the classic Siegel paper \cite{Siegel}. In this paper we  prove the result of Siegel on the mean values of the auxiliary function 
$\Rzeta(s)$. See \cite{A166} for the definition and properties of $\Rzeta(s)$. 

Siegel uses a particular method to obtain his results.  Siegel asserts that 
Riemann tried to get an expression for the zeros of $\Rzeta(s)$. To this end, Riemann
considers the modulus and the argument of $\Rzeta(s)$. Siegel asserts that 
to study the modulus Riemann consider the integral
\begin{equation}
|\Rzeta(\sigma+it)|^2=\int_{0\swarrow1}\int_{0\searrow1}\frac{x^{-\sigma-ti}y^{-\sigma
+ti}e^{\pi i (x^2-y^2)}}{(e^{\pi i x}-e^{-\pi i x})(e^{\pi i y}-e^{-\pi i y})}\,dx\,dy
\end{equation}
where he changes variables, changes the path of integration, and applies the residue theorem. Although these transformations apparently do not give any useful result.

However, Siegel employs this Riemann technique to get some of his results. 
He uses the formula
\begin{equation}
\int_0^{\infty}|\Rzeta(\sigma+it)|^2 e^{-\varepsilon t}\,dt
=
\int_0^{+\infty}e^{-\varepsilon t}\Bigl\{\int_{0\swarrow1}\int_{0\searrow1}\frac{x^{-\sigma-ti}y^{-\sigma
+ti}e^{\pi i (x^2-y^2)}}{(e^{\pi i x}-e^{-\pi i x})(e^{\pi i y}-e^{-\pi i y})}\,dx\,dy
\Bigr\}\,dt
\end{equation}
and apply Riemann transformations and Fubini's Theorem to get 
\begin{equation}\label{firstequation}
\int_0^{\infty}|\Rzeta(\sigma+it)|^2 e^{-\varepsilon t}\,dt\sim
\frac{1}{2\varepsilon}(2\pi\varepsilon)^{\sigma-\frac12}\Gamma(\tfrac12-\sigma).
\end{equation}
I have not been able to reproduce the Siegel proof. He does not give
any other explanation of the procedure.

From \eqref{firstequation} Siegel asserts that for all
$\sigma<\frac12$
\begin{equation}\label{forsigma}
\int_1^{\infty}|\Rzeta(\sigma+it)|^2\Bigl(\frac{t}{2\pi}\Bigr)^\sigma 
e^{-\varepsilon t}\,dt\sim \frac{(2\varepsilon)^{-\frac32}}{1-2\sigma}.
\end{equation}
He does not give a proof of this implication. There is a statement in  
Titchmarsh \cite{T}*{p.~159} which proves the implication \eqref{firstequation} $\Longrightarrow$ \eqref{forsigma} for $\sigma<0$.  The proof in Titchmarsh needs  hypothesis $\sigma<0$.  Siegel only uses \eqref{forsigma} for a value of $\sigma=-5.47559\dots$  so this will be sufficient. However, in Theorem \ref{T:Siegelstatement} we give a different and complete proof of \eqref{forsigma}  for $\sigma<\frac12$.

Theorems \ref{T:main} and \ref{second}  give the mean values of $|\Rzeta(\sigma+it)|^2(t/2\pi)^{\sigma}$ for all $\sigma\in\R$.  Then Theorem \ref{T:meanvalues} gives the main terms of the mean value of $|\Rzeta(\sigma+it)|^2$. 

In \cite{A166}*{eq.~(52)} we saw that $\Rzeta(\frac12+it)=\frac12 e^{-i\vartheta(t)}(Z(t)+iY(t))$, therefore $2|\Rzeta(\frac12+it)|\ge|\zeta(\frac12+it)|$ with equality when $Y(t)=0$. This implies a simple relationship between the mean values of $\Rzeta(\frac12+it)$ and those of $\zeta(\frac12+it)$. 
\begin{figure}[H]
\begin{center}
\includegraphics[width=0.6\hsize]{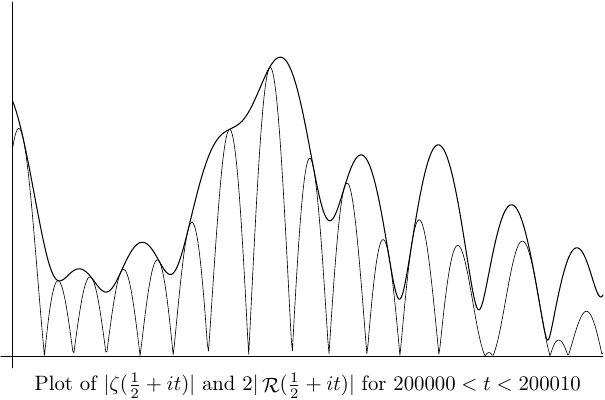}
\end{center}
\end{figure}

\section{Weighted mean values.}

We require the following lemmas.
\begin{lemma}\label{lema1}
For $0<a<b$, and $\alpha$ and $\beta\in\R$ with $\beta\ne0$
\begin{equation}
\Bigl|\int_a^b t^\alpha\cos\beta t\,dt\Bigr|\le \frac{3}{|\beta|}\max(a^\alpha,b^\alpha).
\end{equation}
\end{lemma}

\begin{proof}
We have for $\alpha\ne0$ and $\beta\ne0$
\begin{displaymath}
\int_a^b t^\alpha\cos\beta t\,dt=\frac{b^\alpha\sin(\beta b)-a^\alpha\sin(\beta a)}
{\beta}-\alpha\int_a^b\frac{\sin\beta t}{\beta}t^{\alpha-1}\,dt
\end{displaymath}
and 
\begin{displaymath}
\Bigl|\alpha\int_a^b\frac{\sin\beta t}{\beta}t^{\alpha-1}\,dt\Bigr|\le \frac{|\alpha|}{|\beta|}\int_a^bt^{\alpha-1}\,dt=
\frac{|\alpha|}{\alpha|\beta|}(b^\alpha-a^\alpha)\le \frac{1}{|\beta|}\max(b^\alpha,a^\alpha).
\end{displaymath}
The case $\alpha=0$ is easier.
\end{proof}

\begin{lemma} \label{L:euler}
We have for $x\to+\infty$
\begin{equation}
\sum_{n\le x}\frac{1}{n^{2\sigma}}=\begin{cases}
\frac{x^{1-2\sigma}}{1-2\sigma}+\Orden(x^{-2\sigma}),&\text{for $\sigma\le0$},\\
\noalign{\smallskip}
\zeta(2\sigma)+\frac{x^{1-2\sigma}}{1-2\sigma}+\Orden(x^{-2\sigma}),&\text{for $\sigma>0$, $\sigma\ne1/2$},\\
\noalign{\smallskip}
\log x+\gamma+\Orden(x^{-1}),&\text{for $\sigma=1/2$}.\end{cases}\end{equation}
\end{lemma}
\begin{proof}
This is well known. For example, this can be found in \cite{MV}*{Corollary 1.15}.
\end{proof}

The next two  lemmas are  companions to Lemma 7.2 of Titchmarsh \cite{T}, extending the range of $\sigma$, but without  uniformity in $\sigma$.
\begin{lemma}\label{lema2}
For $\sigma<0$ we have 
\begin{equation}
\sum_{m<n\le x}\frac{n^\sigma}{m^\sigma\log(n/m)}=\Orden(x^2\log x),
\end{equation}
where the implicit constant depend on $\sigma$.
\end{lemma}

\begin{proof}
Let $\Sigma_1$ denote the sum of the terms with $n/m\ge2$ and $\Sigma_2$ the remainder.
In $\Sigma_1$ we have $\log(n/m)\ge\log 2$ and 
\begin{displaymath}
\Sigma_1\le A\sum_{m<n\le x}\frac{n^\sigma }{m^\sigma}\le A2^\sigma 
\sum_{m<n\le x}1=\Orden(x^2).
\end{displaymath}
On $\Sigma_2$ put $n=m+r$ where $1\le r< n/2$, then 
\begin{displaymath}
\log\frac{n}{m}=-\log\frac{n-r}{n}>\frac{r}{n}.
\end{displaymath}
So,
\begin{displaymath}
\Sigma_2\le \sum_{n\le x, r<n/2}\frac{n^\sigma }{(n-r)^\sigma}\frac{n}{r}<
\sum_{n\le x, r<n/2}\frac{n}{r}\ll \sum_{n\le x} n\log n=\Orden(x^2\log x).\qedhere
\end{displaymath}
\end{proof}

\begin{lemma}\label{lema3}
With implicit constants depending on $\sigma$, we have 
\begin{equation}
\sum_{m<n\le x}\frac{1}{n^\sigma m^\sigma \log(n/m)}=
\begin{cases}
\Orden(x^{2-2\sigma}\log x) & \text{if $\sigma<1$,}\\
\Orden(\log^2 x) & \text{if $\sigma=1$,}\\
\Orden(1) & \text{if $\sigma>1$.}
\end{cases}
\end{equation}
\end{lemma}

\begin{proof}  
Divide the sum in $\Sigma_1$ and $\Sigma_2$ as in the proof of the lemma \ref{lema2}. 
For $\Sigma_1$ when $n/m\ge2$  we have
\begin{displaymath}
\Sigma_1\le A \sum_{m<n\le x}\frac{1}{n^{\sigma}m^\sigma}\le 
A\Bigl(\sum_{n\le x}\frac{1}{n^\sigma}\Bigr)^2=
\begin{cases}
\Orden(x^{2-2\sigma}) & \text{for $\sigma<1$,}\\
\Orden(\log ^2x) & \text{for $\sigma=1$,}\\
\Orden(1) & \text{for $\sigma>1$.}
\end{cases}
\end{displaymath}
We apply to $\Sigma_2$ the same procedure as in Lemma \ref{lema2} (with $n/2<(n-r)<n$)
\[\Sigma_2\le \sum_{n\le x, r<n/2}\frac{n}{n^\sigma(n-r)^\sigma r}\ll  \sum_{n\le x} n^{1-2\sigma}\log n=
\begin{cases}
\Orden(x^{2-2\sigma}\log x) & \text{for $\sigma<1$,}\\
\Orden(\log ^2x) & \text{for $\sigma=1$,}\\
\Orden(1) & \text{for $\sigma>1$.}
\end{cases}
\]
Recall that our constants depends on $\sigma$. 
\end{proof}

\begin{theorem}\label{T:main}
When $T\to+\infty$, the integral
\begin{equation}
\frac{1}{T}\int_1^T|\Rzeta(\sigma+it)|^2\Bigl(\frac{t}{2\pi}\Bigr)^{\sigma}\,dt
\end{equation}
with $\sigma\ge0$ is equal to 
\begin{equation}
\begin{cases}
\frac{2}{3(1-2\sigma)}\left(\frac{T}{2\pi}\right)^{1/2}+
\Orden(T^\frac14) & \text{if $0\le \sigma\le \tfrac14$,}\\ \noalign{\smallskip}
\frac{2}{3(1-2\sigma)}\left(\frac{T}{2\pi}\right)^{1/2}+
\frac{\zeta(2\sigma)}{\sigma+1}\left(\frac{T}{2\pi}\right)^{\sigma}+
\Orden(T^\frac14) & \text{if $1/4< \sigma<\tfrac12$,}\\ \noalign{\smallskip}
\frac{1}{3}\left(\frac{T}{2\pi}\right)^{1/2}\log\frac{T}{2\pi}+
2\frac{3\gamma-4}{9}\left(\frac{T}{2\pi}\right)^{1/2}+\Orden(T^{\frac14}\sqrt{\log T})
 & \text{if $\sigma=1/2$,}\\\noalign{\smallskip}
\frac{\zeta(2\sigma)}{\sigma+1}\left(\frac{T}{2\pi}\right)^{\sigma}+
\frac{2}{3(1-2\sigma)}\left(\frac{T}{2\pi}\right)^{1/2}+
\Orden(T^{\frac{\sigma}{2}}) & \text{if $1/2<\sigma<1$,}\\ \noalign{\smallskip}
\frac{\zeta(2\sigma)}{\sigma+1}\left(\frac{T}{2\pi}\right)^{\sigma}+
\Orden(T^{\frac{\sigma}{2}}) & \text{if $1\le\sigma\le 2$,}\\ \noalign{\smallskip}
\frac{\zeta(2\sigma)}{\sigma+1}\left(\frac{T}{2\pi}\right)^{\sigma}+
\Orden(T^{\sigma-1}) & \text{if $2<\sigma$.}
\end{cases}
\end{equation}
where the implicit constants depend on $\sigma$ but of course, not on $T$. 
\end{theorem}

\begin{proof}
We start from the result proved in  \cite{A86}. For any value of $\sigma\in\R$ and  $t>0$ we have for $t\to+\infty$
\begin{equation}\label{main}
\Rzeta(\sigma+it)=\sum_{n\le \sqrt{t/2\pi}}\frac{1}{n^{\sigma+it}}+\Orden(t^{-\sigma/2}).
\end{equation}
We call $S(t)$ the first sum and $R(t)$ the term of error. 
We have 
\begin{equation}\label{normR}
\Vert R\Vert_2:=\Bigl(\frac{1}{T}\int_1^T|R(t)|^2(t/2\pi)^\sigma\,dt\Bigr)^{1/2}=
\Orden(1).
\end{equation}
(All implicit  constants in the $\Orden$ terms may depend on $\sigma$.)

And we have
\begin{multline}\label{firststep}
\int_1^T|S(t)|^2\left(\frac{t}{2\pi}\right)^{\sigma}\,dt=\sum_{n\le \sqrt{T/2\pi}}
\frac{1}{n^{2\sigma}}\int_{2\pi n^2}^T\left(\frac{t}{2\pi}\right)^{\sigma}\,dt+\\
+
2\sum_{m<n\le \sqrt{T/2\pi}}\frac{1}{n^\sigma m^\sigma}\int_{2\pi n^2}^T
\left(\frac{t}{2\pi}\right)^{\sigma}\cos(t\log(n/m))\,dt:=S_1+S_2
\end{multline}
say. Applying the lemma \ref{lema1} to the second sum $S_2$ in \eqref{firststep} yields
\begin{displaymath}
|S_2|\le\frac{6T^\sigma}{(2\pi)^\sigma}\sum_{m<n\le\sqrt{T/2\pi}}\frac{1}{n^\sigma m^\sigma\log(n/m)}
\end{displaymath}
Now we apply Lemma \ref{lema3} to get
\begin{equation}
|S_2|=\Orden(T\log T), \quad \Orden(T\log^2T),\quad \Orden(T^\sigma),
\end{equation}
respectively for $0\le \sigma<1$, $\sigma=1$ or $\sigma>1$. 
The integral in the first sum may be computed explicitly. Since we assume $\sigma\ge0$ we have
\begin{multline}\label{integrated}
S_1=(2\pi)^{-\sigma}\sum_{n\le\sqrt{T/2\pi}}\frac{T^{\sigma+1}-(2\pi
n^2)^{\sigma+1}}{n^{2\sigma}(\sigma+1)}=\\ =
\frac{T}{\sigma+1}\Bigl(\frac{T}{2\pi}\Bigr)^{\sigma}\sum_{n\le\sqrt{T/2\pi}}
\frac{1}{n^{2\sigma}}-\frac{2\pi}{\sigma+1}\sum_{n\le\sqrt{T/2\pi}}n^2.
\end{multline}
For  $0\le\sigma<1/2$ we have  (by Lemma \ref{L:euler})
\[
S_1=\frac{T}{\sigma+1}\Bigl(\frac{T}{2\pi}\Bigr)^{\sigma}
\Bigl(\frac{(T/2\pi)^{(1-2\sigma)/2}}{1-2\sigma}+\zeta(2\sigma)+
\Orden(T^{-\sigma})\Bigr)
-\frac{2\pi}{\sigma+1}\Bigl(
\frac{(T/2\pi)^{3/2}}{3}+\Orden(T)\Bigr).
\]
Rearranging the terms, we obtain
\begin{displaymath}
S_1=\frac{2T}{3(1-2\sigma)}\Bigl(\frac{T}{2\pi}\Bigr)^{\tfrac12}+
\frac{\zeta(2\sigma)T}{\sigma+1}\Bigl(\frac{T}{2\pi}\Bigr)^{\sigma}+\Orden(T).
\end{displaymath}
For $\sigma=\tfrac12$,
\begin{displaymath}
S_1=\frac{2T}{3}\Bigl(\frac{T}{2\pi}\Bigr)^{\frac12}
\Bigl(\frac{1}{2}\log\frac{T}{2\pi}+\gamma+\Orden(T^{-1/2})
\Bigr)-\frac{4\pi}{3}\Bigl(
\frac{(T/2\pi)^{3/2}}{3}+\Orden(T)\Bigr).
\end{displaymath}
Ordering the terms yields
\begin{displaymath}
S_1=\frac{T
}{3}\Bigl(\frac{T}{2\pi}\Bigr)^{\frac12}\log\frac{T}{2\pi}+2T\frac{3\gamma-1}{9}
\Bigl(\frac{T}{2\pi}\Bigr)^{\frac12}+\Orden(T).
\end{displaymath}
For $\sigma>\tfrac12$
\[
S_1=
\frac{T}{\sigma+1}\Bigl(\frac{T}{2\pi}\Bigr)^{\sigma}\Bigl(\zeta(2\sigma)+\frac{1}{1-2\sigma}
\Bigl(\frac{T}{2\pi}\Bigr)^{\frac12-\sigma}+\Orden(T^{-\sigma})\Bigr) -
\frac{2\pi}{\sigma+1}\Bigl(\frac{(T/2\pi)^{3/2}}{3}+\Orden(T)\Bigr).
\]
Rearranging the terms, we get
\begin{equation}
S_1=\frac{T\zeta(2\sigma)}{\sigma+1}\Bigl(\frac{T}{2\pi}\Bigr)^{\sigma}+\frac{2T}{3(1-2\sigma)}\Bigl(\frac{T}{2\pi}\Bigr)^{\frac12}+\Orden(T).
\end{equation}
From the above, it follows that 
\begin{displaymath}
\Vert S\Vert_2=\Bigl(\frac1T\int_1^T|S(t)|^2(t/2\pi)^\sigma\,dt
\Bigr)^{1/2}=
\begin{cases}
\Orden(T^{1/4}) &\text{for $0\le \sigma< \tfrac12$,}\\
\Orden(T^{1/4}\sqrt{\log T}) &\text{for $\sigma=\tfrac12$,}\\
\Orden(T^{\sigma/2}) &\text{for $\tfrac12<\sigma$-}
\end{cases}
\end{displaymath}
Therefore,
\begin{equation}\label{scalarproduct}
\frac1T\int_1^T|S(t)R(t)|(t/2\pi)^\sigma\,dt
=
\begin{cases}
\Orden(T^{1/4}) &\text{for $0\le \sigma< \tfrac12$,}\\
\Orden(T^{1/4}\sqrt{\log T}) &\text{for $\sigma=\tfrac12$,}\\
\Orden(T^{\sigma/2}) &\text{for $\tfrac12<\sigma$.}
\end{cases}
\end{equation}
By \eqref{normR} and \eqref{scalarproduct}   we have 
\begin{multline*}
\frac{1}{T}\int_1^T|\Rzeta(\sigma+it)|^2\Bigl(\frac{t}{2\pi}\Bigr)^{\sigma}\,dt=
\Vert \Rzeta\Vert_2^2=\Vert S\Vert_2^2+\Vert R\Vert_2^2+2\Re\langle R,S\rangle
=\\ =
\begin{cases}
\Vert S\Vert_2^2+\Orden(T^{\frac14})&\text{for $0\le\sigma<\frac12$,}\\
\Vert S\Vert_2^2+\Orden(T^{\frac14}\sqrt{\log T})&\text{for $\sigma=\tfrac12$,}\\
\Vert S\Vert_2^2+\Orden(T^{\frac{\sigma}{2}})&\text{for $\tfrac12<\sigma$.}
\end{cases}
\end{multline*}
Now 
\begin{displaymath}
\Vert S\Vert_2^2\le \frac{1}{T}(S_1+S_2)=\frac{S_1}{T}+\Orden(\log T),\quad 
\frac{S_1}{T}+\Orden(\log^2 T),\quad \frac{S_1}{T}+\Orden(T^{\sigma-1}),
\end{displaymath}
respectively in the three cases $0\le \sigma<1$, $\sigma=1$ or $\sigma>1$.
Joining this with the values of $S_1$ given above, we get for $0\le \sigma<\tfrac12$
\begin{displaymath}
\frac{1}{T}\int_1^T|\Rzeta(\sigma+it)|^2\Bigl(\frac{t}{2\pi}\Bigr)^{\sigma}\,dt=
\frac{2}{3(1-2\sigma)}\Bigl(\frac{T}{2\pi}\Bigr)^{\frac12}+
\frac{\zeta(2\sigma)}{\sigma+1}\Bigl(\frac{T}{2\pi}\Bigr)^{\sigma}+\Orden(T^{\frac14}).
\end{displaymath}
For $\sigma=\tfrac12$ we have
\begin{displaymath}
\frac{1}{T}\int_1^T|\Rzeta(\sigma+it)|^2\Bigl(\frac{t}{2\pi}\Bigr)^{\frac12}\,dt
=
\frac{1
}{3}\Bigl(\frac{T}{2\pi}\Bigr)^{\frac12}\log\frac{T}{2\pi}+2\frac{3\gamma-1}{9}
\Bigl(\frac{T}{2\pi}\Bigr)^{\frac12}+\Orden(T^{\frac14}\sqrt{\log T}).
\end{displaymath}
For $\tfrac12<\sigma<1$ we have 
\begin{displaymath}
\frac{1}{T}\int_1^T|\Rzeta(\sigma+it)|^2\Bigl(\frac{t}{2\pi}\Bigr)^{\sigma}\,dt=
\frac{\zeta(2\sigma)}{\sigma+1}\Bigl(\frac{T}{2\pi}\Bigr)^{\sigma}+
\frac{2}{3(1-2\sigma)}\Bigl(\frac{T}{2\pi}\Bigr)^{\frac12}+
\Orden(T^{\frac{\sigma}{2}});
\end{displaymath}
for $1\le \sigma\le 2$, the last term is less than the error and we only get
\begin{displaymath}
\frac{1}{T}\int_1^T|\Rzeta(\sigma+it)|^2\Bigl(\frac{t}{2\pi}\Bigr)^{\sigma}\,dt=
\frac{\zeta(2\sigma)}{\sigma+1}\Bigl(\frac{T}{2\pi}\Bigr)^{\sigma}+
\Orden(T^{\frac{\sigma}{2}})
\end{displaymath}
For $\sigma>2$ the error from $S_2/T$ dominates, and we have 
\begin{displaymath}
\frac{1}{T}\int_1^T|\Rzeta(\sigma+it)|^2\Bigl(\frac{t}{2\pi}\Bigr)^{\sigma}\,dt=
\frac{\zeta(2\sigma)}{\sigma+1}\Bigl(\frac{T}{2\pi}\Bigr)^{\sigma}+
\Orden(T^{\sigma-1}).\qedhere
\end{displaymath}
\end{proof}

\begin{theorem}\label{second}
For $\sigma<0$  we have
\begin{equation}
\frac{1}{T}\int_1^T|\Rzeta(\sigma+it)|^2\Bigl(\frac{t}{2\pi}\Bigr)^{\sigma}\,dt=
\frac{2}{3(1-2\sigma)}\Bigl(\frac{T}{2\pi}\Bigr)^{\frac12}+\Orden(T^{\frac14}).
\end{equation}
\end{theorem}

\begin{proof}
Again, we start from \eqref{main} and call the two terms $S$ and $R$. 
As before, we have
\begin{displaymath}
\Vert R\Vert_2=\Bigl(\frac{1}{T}
\int_1^T|\Rzeta(\sigma+it)|^2\Bigl(\frac{t}{2\pi}\Bigr)^{\sigma}\,dt\Bigr)^{\frac12}
=\Orden(1).
\end{displaymath}

Now, the integral for $|S(t)|^2$ is given as before by  \eqref{firststep}. But now by Lemma \ref{lema1} the sum
$S_2$ is bounded by
\begin{displaymath}
|S_2|\le \frac{6}{(2\pi)^\sigma}\sum_{m<n\le \sqrt{T/2\pi}}\frac{(2\pi n^2)^\sigma}{n^\sigma m^\sigma\log (n/m)},
\end{displaymath}
so that by Lemma \ref{lema2}
\begin{equation}
|S_2|=\Orden(T\log T).
\end{equation}
We assume $\sigma<0$ when $\sigma\ne -1$ we have as before \eqref{integrated}, 
so that for $\sigma\ne-1$ and $\sigma<0$
\[
S_1=\frac{T}{\sigma+1}\Bigl(\frac{T}{2\pi}\Bigr)^{\sigma}
\Bigl(\frac{(T/2\pi)^{(1-2\sigma)/2}}{1-2\sigma}+
\Orden(T^{-\sigma})\Bigr)
-\frac{2\pi}{\sigma+1}\Bigl(
\frac{(T/2\pi)^{3/2}}{3}+\Orden(T)\Bigr).
\]
Rearranging the terms, this is 
\begin{equation}\label{S1result}
S_1=\frac{2T}{3(1-2\sigma)}\Bigl(\frac{T}{2\pi}\Bigr)^{\frac12}+\Orden(T).
\end{equation}

When $\sigma=-1$, the integral in \eqref{firststep} with $x=\sqrt{T/2\pi}$ is 
\begin{displaymath}
S_1=2\pi\sum_{n\le \sqrt{T/2\pi}}n^2\log\frac{T}{2\pi n^2}=4\pi\sum_{n\le x}
n^2\log \frac{x}{n},
\end{displaymath}
integrating by parts transforms into
\begin{displaymath}
S_1=\frac{4\pi}{3}\int_1^x(t^3+\Orden(t^2))\frac{dt}{t}=\frac{4\pi}{9}x^3+\Orden(x^2).
\end{displaymath}
It follows that 
\begin{displaymath}
S_1=\frac{2T}{9}\Bigl(\frac{T}{2\pi}\Bigr)^{\frac12}+\Orden(T).
\end{displaymath}
Thus, \eqref{S1result} is also valid for $\sigma=-1$. 

Then 
\begin{displaymath}
\Vert S\Vert_2=\Bigl(\frac1T\int_1^T|S(t)|^2(t/2\pi)^\sigma\,dt
\Bigr)^{1/2}=\Orden(T^{\frac14}),
\end{displaymath}
and
\begin{displaymath}
\frac1T\int_1^T|S(t)R(t)|(t/2\pi)^\sigma\,dt\le \Vert S\Vert_2\,\Vert R\Vert_2
\le \Orden(T^{\frac14}).
\end{displaymath}
So,
\begin{align*}
\frac{1}{T}\int_1^T|\Rzeta(\sigma+it)|^2\Bigl(\frac{t}{2\pi}\Bigr)^{\sigma}\,dt&=
\Vert \Rzeta\Vert_2^2= \Vert S\Vert_2^2+\Vert R\Vert_2^2+2\Re\langle R,S\rangle
\\&
=\frac{S_1+S_2}{T}+\Orden(T^{\frac14}) =
\frac{2}{3(1-2\sigma)}\Bigl(\frac{T}{2\pi}\Bigr)^{\frac12}+\Orden(T^{\frac14}).\qedhere
\end{align*}
\end{proof}

\begin{theorem}\label{T:Siegelstatement}
For all $\sigma<\tfrac12$ we have
\begin{equation}\label{forsigma2}
\int_1^{\infty}|\Rzeta(\sigma+it)|^2\Bigl(\frac{t}{2\pi}\Bigr)^\sigma 
e^{-\varepsilon t}\,dt= \frac{(2\varepsilon)^{-\frac32}}{1-2\sigma}+\orden(\varepsilon^{-3/2}).
\end{equation}
\end{theorem}

\begin{proof}
From Theorems \ref{T:main} and \ref{second} we get that for all $\sigma<\tfrac12$ 
\begin{displaymath}
F(T):=\int_1^T|\Rzeta(\sigma+it)|^2\Bigl(\frac{t}{2\pi}\Bigr)^{\sigma}\,dt=
\frac{2}{3(1-2\sigma)\sqrt{2\pi}}T^{\frac32}+\Orden(T^a),
\end{displaymath}
where $a=\max(\frac54,1+\sigma)<\frac32$.

Then the integral in \eqref{forsigma2} is equal to 
\begin{multline*}
\int_1^\infty e^{-\varepsilon t}\,dF(t)=\varepsilon
\int_1^\infty F(t)e^{-\varepsilon t}\,dt=\\
=\frac{2\varepsilon}{3(1-2\sigma)\sqrt{2\pi}}\int_1^\infty t^{3/2}e^{-\varepsilon t}\,dt+
\Orden\Bigl(\frac{2\varepsilon}{3(1-2\sigma)\sqrt{2\pi}}\int_1^\infty t^{a}e^{-\varepsilon t}\,dt
\Bigr)
=\frac{(2\varepsilon)^{-\frac32}}{1-2\sigma}+\Orden(\varepsilon^{-a}).
\end{multline*}
Since it is easy to see that for $\varepsilon>0$ and $a>0$
\[\int_1^\infty t^a e^{-\varepsilon t}\,dt=\frac{\Gamma(1+a)}{\varepsilon^{1+a}}+\sum_{n=1}^\infty (-1)^n \frac{\varepsilon^{n-1}}{(n-1)!\;(a+n)}.\qedhere\]
\end{proof}

\section{Mean values}

\begin{theorem}\label{T:meanvalues}
When $T\to+\infty$, we have 
\begin{equation}\label{E:meanvalues}
\frac{1}{T}\int_0^T|\Rzeta(\sigma+it)|^2\,dt
=\begin{cases}
\frac{2}{(1-2\sigma)(3-2\sigma)}\left(\frac{T}{2\pi}\right)^{\frac12-\sigma}+\Orden(T^{\frac14-\sigma}),&\text{for $\sigma\le 1/4$},\\ \noalign{\smallskip}
\frac{2}{(1-2\sigma)(3-2\sigma)}\left(\frac{T}{2\pi}\right)^{\frac12-\sigma}+\zeta(2\sigma)+\Orden(T^{\frac14-\sigma}),&\text{for $1/4<\sigma< 1/2$},\\ \noalign{\smallskip}
\frac{1}{2}\log\frac{T}{2\pi}+(\gamma-\tfrac12)+\Orden(T^{-1/4}\log^{\frac12}T),&\text{for $\sigma=1/2$},\\ \noalign{\smallskip}
\zeta(2\sigma)+\frac{2}{(1-2\sigma)(3-2\sigma)}\left(\frac{T}{2\pi}\right)^{\frac12-\sigma}+\Orden(T^{-\frac{\sigma}{2}}),&\text{for $1/2<\sigma< 1$},\\ \noalign{\smallskip}
\zeta(2)+\Orden(T^{-\frac{1}{2}}\log^{\frac12}T),&\text{for $\sigma=1$},\\ \noalign{\smallskip}
\zeta(2\sigma)+\Orden(T^{-\frac{\sigma}{2}}),&\text{for $1<\sigma< 2$},\\ \noalign{\smallskip}
\zeta(2\sigma)+\Orden(T^{-1}),&\text{for $2\le \sigma$}.
\end{cases}
\end{equation}
\end{theorem}

\begin{proof}
As in \eqref{main} we have $\Rzeta(\sigma+it)=S(t)+R(t)$ with 
\begin{equation}\label{normR2}
\Vert R\Vert_2:=\Bigl(\frac{1}{T}\int_0^T|R(t)|^2\,dt\Bigr)^{1/2}=
\begin{cases}
\Orden(T^{-\sigma/2}), & \text{for $\sigma\ne1$},\\
\Orden(T^{-1/2}\log^{1/2} T), & \text{for $\sigma=1$}.\end{cases}
\end{equation}
The term in $S(t)$ gives us
\begin{multline}\label{firststep2}
\int_1^T|S(t)|^2\,dt=\sum_{n\le \sqrt{T/2\pi}}
\frac{1}{n^{2\sigma}}\int_{2\pi n^2}^T\,dt+\\
+
2\sum_{m<n\le \sqrt{T/2\pi}}\frac{1}{n^\sigma m^\sigma}\int_{2\pi n^2}^T
\cos(t\log(n/m))\,dt:=S_1+S_2
\end{multline}
say.  Applying the lemma \ref{L:euler} we get 
\[
S_1=\begin{cases}
\frac{2T}{(1-2\sigma)(3-2\sigma)}\left(\frac{T}{2\pi}\right)^{\frac12-\sigma}+\Orden(T^{1-\sigma}),&\text{for $\sigma\le 0$},\\ \noalign{\smallskip}
\frac{2T}{(1-2\sigma)(3-2\sigma)}\left(\frac{T}{2\pi}\right)^{\frac12-\sigma}+\zeta(2\sigma)T+\Orden(T^{1-\sigma}),&\text{for $0<\sigma< 1/2$},\\ \noalign{\smallskip}
\frac{T}{2}\log\frac{T}{2\pi}+(\gamma-\tfrac12)T+\Orden(T^{1/2}),&\text{for $\sigma=1/2$},\\ \noalign{\smallskip}
\zeta(2\sigma)T+\frac{2T}{(1-2\sigma)(3-2\sigma)}\left(\frac{T}{2\pi}\right)^{\frac12-\sigma}+\Orden(T^{1-\sigma}),&\text{for $1/2<\sigma\le 1$},\\ \noalign{\smallskip}
\zeta(2\sigma)T+\frac{2T}{(1-2\sigma)(3-2\sigma)}\left(\frac{T}{2\pi}\right)^{\frac12-\sigma}+\Orden(T^{1-\sigma}),&\text{for $1<\sigma< 3/2$},\\ \noalign{\smallskip}
\zeta(3)T-\pi\log\frac{T}{2\pi}-2\pi(\gamma+\tfrac12)+\Orden(T^{-1/2}),&\text{for $\sigma=3/2$},\\ \noalign{\smallskip}
\zeta(2\sigma)T-2\pi\zeta(2\sigma-2)+\frac{2T}{(1-2\sigma)(3-2\sigma)}\left(\frac{T}{2\pi}\right)^{\frac12-\sigma}+\Orden(T^{1-\sigma}),&\text{for $\sigma\ge 3/2$}.
\end{cases}
\]
Applying Lemma \ref{lema1} to $S_2$, and then by Lemma \ref{lema3} we find
\[|S_2|\le \sum_{m<n\le \sqrt{T/2\pi}}\frac{6}{n^\sigma m^\sigma\log(n/m)}\le \begin{cases}
\Orden (T^{1-\sigma}\log T), &\text{for $\sigma<1$},\\
\Orden (\log^2 T), &\text{for $\sigma=1$},\\
\Orden (1), &\text{for $1<\sigma$}.\\
\end{cases}\]

Then 
\begin{displaymath}
\Vert S\Vert_2=\Bigl(\frac1T\int_1^T|S(t)|^2\,dt
\Bigr)^{1/2}=\Orden(T^{\frac14-\frac{\sigma}{2}}),\quad \Orden(\sqrt{\log T}),\quad \Orden(1), 
\end{displaymath}
for $\sigma<\frac12$, $\sigma=\frac12$, and $\sigma>\frac12$ respectively.
Therefore,  
\begin{multline*}
\frac1T\int_1^T|S(t)R(t)|\,dt\le \Vert S\Vert_2\,\Vert R\Vert_2\\
\le \Orden(T^{\frac14-\sigma}),\quad \Orden(T^{-\frac{\sigma}{2}}\log^{1/2} T),\quad \Orden(T^{-\frac{\sigma}{2}}), \quad \Orden(T^{-1/2}\log^{1/2}T),\quad \Orden(T^{-\frac{\sigma}{2}}),
\end{multline*}
for $\sigma<1/2$, $\sigma=\frac12$, $\frac12<\sigma<1$, $\sigma=1$   and $\sigma>1$ respectively.
So,
\begin{multline*}
\frac{1}{T}\int_1^T|\Rzeta(\sigma+it)|^2\,dt=
\Vert \Rzeta\Vert_2^2= \Vert S\Vert_2^2+\Vert R\Vert_2^2+2\Re\langle R,S\rangle
=\\
=\frac{S_1+S_2}{T}+\Vert R\Vert_2^2+2\Re\langle R,S\rangle.
\end{multline*}
The above results allow us to bound 
\[\frac{S_2}{T}+\Vert R\Vert_2^2+\Vert R\Vert_2\Vert S\Vert_2=
\begin{cases}
\Orden(T^{\frac14-\sigma}),&\text{for $\sigma<1/2$},\\
\Orden(T^{-\frac14}\log^{1/2}T),&\text{for $\sigma=1/2$},\\
\Orden(T^{-\sigma/2}),&\text{for $1/2<\sigma<1$},\\
\Orden(T^{-1/2}\log^{1/2}T),&\text{for $\sigma=1$},\\
\Orden(T^{-\sigma/2}),&\text{for $1<\sigma<2$},\\
\Orden(T^{-1}),&\text{for $2\le \sigma$}.
\end{cases}
\]
Combining with the values of $S_1/T$ obtained above yields \eqref{E:meanvalues}.
\end{proof}

\end{document}